\documentclass[11pt]{article}
\usepackage{amsfonts,mathrsfs}
\usepackage{multicol}
\usepackage{geometry}
\usepackage{float}
\usepackage{tikz,chemfig}
\usepackage{yhmath}
\usepackage{indentfirst}   
\usepackage{epic}
\usepackage{etex}
\usepackage{pinlabel}
\usepackage{ulem}
\usepackage{cases}            
\usepackage{setspace,color}
\usepackage{enumerate,mathdots}
\usepackage{amsmath,amssymb,bm} 
\usepackage{graphics}
\usepackage{epsf, epstopdf}
\usepackage{color, ifpdf}  
\usepackage{graphicx,booktabs,multirow}
\usepackage{latexsym}
\usepackage{appendix}
\usepackage{hyperref}
\usepackage{skak}
\usepackage{authblk}
\usepackage[font=small,labelfont=bf]{caption}

\newcommand\red[1] {{\color{red} #1}}

\newcommand\blue[1] {{\color{blue} #1}}

\newcommand\magenta[1] {{\color{magenta} #1}}

\usetikzlibrary{backgrounds}
\usetikzlibrary{topaths,calc}
\usetikzlibrary{arrows}
\usetikzlibrary{shapes,shapes.geometric,shapes.misc}
\pgfdeclarelayer{edgelayer}
\pgfdeclarelayer{nodelayer}
\pgfsetlayers{background,edgelayer,nodelayer,main}
\tikzstyle{none}=[inner sep=0mm]

\newcommand \tikzp[2]
{
	\begin{center}
		\begin{tikzpicture}[scale=#1]
			#2
		\end{tikzpicture}
	\end{center}
}

\tikzstyle{cblue}=[circle, draw, thin,fill=blue!20, scale=0.5]

\newcommand \equ[2]
{
	\begin{equation}\label{#1}
		#2
	\end{equation}
}

\newtheorem{theo}{Theorem}
\newtheorem{lem}{Lemma}
\newtheorem{pro}{Proposition}
\newtheorem{cor}{Corollary}
\newtheorem{obs}{Observation}
\newtheorem{deff}{Definition}
\newtheorem{prob}{Problem}
\newtheorem{con}{Conjecture}
\newtheorem{clm}{Claim}

\newcommand \them[2]
{\begin{theo}
		\label{#1} #2
	\end{theo}
}

\newcommand \lemm[2]
{\begin{lem}
		\label{#1} #2
	\end{lem}
}

\newcommand \obsv[2]
{\begin{obs}
		\label{#1} #2
	\end{obs}
}

\newcommand \defn[2]
{\begin{deff}
		\label{#1} #2
	\end{deff}
}

\newcommand \corr[2]
{\begin{cor}
		\label{#1} #2
	\end{cor}
}

\newcommand \Clm[2]
{
	\begin{clm}\label{#1}
		#2
	\end{clm}
}

\newcommand \low[2]
{\ell(#1,#2)}

\newcounter{countcase}
\def\incase{\addtocounter{countcase}{1}
	{\noindent {\bf Case \thecountcase}: }}

\newcounter{countclaim}

\def \proof {\noindent {{\it Proof}}.
}

\newcommand{\claimend}{{\hfill $\natural $}}

\newcommand{\proofend}{{\hfill$\Box$}\setcounter{countclaim}{0}\setcounter{countcase}{0}}

\def \N {{\mathbb N}}

\def \hyf {{\cal F}}
\def \hyp {{\cal P}}
\def \hyh {{\cal H}}
\def \hyv {{\cal V}}
\def \hyc {{\cal C}}
\def \hyd {{\cal D}}
\def \hye {{\cal E}}
\def \hyk {{\cal K}}

\newcommand \vtx[2]
{v_{#1,#2}}

\begin{document}
	
\title{
		On the colorability of bi-hypergraphs
	}
	
\author[1]{\small Meiqiao Zhang\thanks{Email: nie21.zm@e.ntu.edu.sg and 
		meiqiaozhang95@163.com}}
\author[1]{\small Fengming Dong\thanks{Corresponding Author. Email: fengming.dong@nie.edu.sg and donggraph@163.com}}
\author[2]{\small Ruixue Zhang\thanks{Email: ruixuezhang7@163.com}}
	
	\affil[1]{\footnotesize 
		National Institute of Education,
	Nanyang Technological University, 
Singapore}
	\affil[2]{\footnotesize School of Mathematics and Statistics, Qingdao University, China}

	\date{}
	
	\maketitle{}

	\abstract{
	A {\it mixed hypergraph}
	$\hyh=(\hyv,\hyc,\hyd)$ consists of the vertex set $\hyv$ and 
	two families of subsets of $2^{\hyv}$:
	the family $\hyc$ of co-edges and 
	the family $\hyd$ of edges.
	$\hyh$ is said to be {\it colorable}
	if there is a mapping $f$ from $\hyv$
	to the set  of positive integers
	such that $|\{f(v):v\in e\}|<|e|$ for each $e\in \hyc$ and $|\{f(v):v\in e\}|>1$ for each $e\in \hyd$.
	There exist mixed hypergraphs 
	which are uncolorable, and quite little about these mixed hypergraphs is known. 
	
	A mixed hypergraph is called 
	a {\it bi-hypergraph} if its co-edge set 
	and edge set are the same.
	In this article, we first apply 
	Lov\'asz local lemma to show that 
	any $r$-uniform bi-hypergraph with $r\ge 4$
	is colorable if every edge is incident to
	less than $(r-1)^{r-1}e^{-1}-1$ other edges, where $e$ is 
	 the base of natural logarithms.
Then, we show that among all the
uncolorable  $3$-uniform bi-hypergraphs, the smallest size of a minimal 
 one is ten, which answers a question raised by Tuza and Voloshin in 2000. As an extension, we provide a minimal uncolorable $3$-uniform bi-hypergraph of order $n$ and size at most $\frac{7n}3-4$ for every $n\ge 6$.
	}


\section{Introduction
\label{intro}}

A hypergraph $\hyh$ consists of a vertex set and an edge set, where any member in the edge set is a subset of its vertex set.
As a natural generalization of the concept of the chromatic number of a graph, Erd\H{o}s and Hajnal \cite{Erdos1966} introduced the chromatic number $\chi(\hyh)$
of a hypergraph $\hyh$
as the minimal number of colors needed to color the vertices in $\hyh$ 
such that no edge in $\hyh$ has 
all its vertices with the same color.

Voloshin \cite{volo1995,volo1993}
extended the concept of the chromatic number of a hypergraph to the upper chromatic number of a mixed hypergraph.
A {\it mixed hypergraph} is a triple
$\hyh=(\hyv, \hyc, \hyd)$,
where $\hyv$ is a finite set and 
both $\hyc$ and $\hyd$ 
are subsets of the power set $2^{\hyv}$.  
The members of $\hyv$ are called 
{\it vertices} of $\hyh$, 
and the members in  $\hyc$ 
and  the members in  $\hyd$
are called 
 {\it co-edges} and {\it edges}
 of $\hyh$ respectively. 
 In this paper, we always assume 
 that  
 both $\hyc$  and   $\hyd$
 are {\it Sperner families}
 (i.e., $e_1\not\subseteq e_2$ for any 
 distinct $e_1,e_2\in \hyc$
 or $e_1,e_2\in \hyd$). 
 
Let $\N$ denote the set of positive integers, and for any $k\in \N$, 
 	let $[k]$ denote the set 
 	$\{1,2,\dots,k\}$.

 \begin{deff}\label{def0} 
Let $\hyh=(\hyv, \hyc, \hyd)$ be a mixed hypergraph and $f$ be a mapping from $\hyv$ to $\N$. For any $e$ in $\hyc$ (or resp., $\hyd$), $e$ is properly colored in $f$ if $|f(e)|<|e|$ (or resp., $|f(e)|>1$), where $f(e)$ denotes the set $\{f(v):v\in e\}$. Then $f$ is a proper coloring of $\hyh$ if every member of $\hyc\cup\hyd$ is properly colored in $f$, and $\hyh$ is colorable if $\hyh$ has a proper coloring. In this case,
 	the maximum number of colors that can occur simultaneously in a
 	proper coloring of $\hyh$ is called 
 	the {\it upper chromatic number}
 	$\bar {\chi}(\hyh)$ of $\hyh$.  
 \end{deff}

 If $\hyc=\emptyset$, then $\hyh=(\hyv, \hyc,\hyd)$ is a hypergraph with edge set $\hyd$. 
 Thus, a proper coloring of a mixed hypergraph is an extension of 
 a proper coloring of a hypergraph. 
 If $\hyd=\emptyset$ and $\hyh$ is uniform,
 $\bar {\chi}(\hyh)$ is treated as a particular case of a more general problem raised by Ahlswede
 et al.~\cite{ahls1992},
 where the authors studied 
 the largest integer, denoted by $N(\hyh,k)$, 
 such that $\hyh$ has a mapping 
$f:\hyv\mapsto \N$ 
 satisfying the conditions 
 $|f(\hyv)|=|N(\hyh,k)|$
 and $|f(e)|\le k$ for each edge $e\in \hyc$.

A proper coloring of a mixed hypergraph $\hyh$ with a certain number of colors
is also called a {\it strict coloring} of $\hyh$.
In \cite{tuza1998},
the concepts of strict coloring and upper chromatic number are applied
to study the Steiner triple systems (STSs).
 In particular, the authors studied 
 two different kinds of
colorings for STSs. 
In one of them, they viewed 
all blocks of the triple system as co-edges;
such systems will be termed CSTS, referring to the expression `co-hypergraph'
and to indicate that monochromatic blocks are allowed. 
In the other type,  all blocks
are treated as edges and co-edges at the same time, 
called Bi-Steiner triple systems. Moreover, proper colorings of mixed hypergraphs also have wide applications in the studies of resource allocation, parallel computing and DNA molecules, as explained in~\cite{volo2002}.

Unfortunately, there exist mixed hypergraphs which do not admit 
proper colorings, 
i.e., \textit{uncolorable}.
Some obvious examples occur when edges of size two arise. For example, $\hyh$ is uncolorable if there is an co-edge $e$ in $\hyh$ such that every pair of vertices in $e$ form an edge in $\hyh$.
A less trivial class of uncolorable mixed hypergraphs was identified by Tuza and Voloshin~\cite{tuza2000}.
For any $n,l,m\in\N$, let $\hyk(n,l,m)$ be the mixed hypergraph $([n],\binom{[n]}{l},\binom{[n]}{m})$. 
Tuza and Voloshin~\cite{tuza2000} have shown that $\hyk(n,l,m)$ is uncolorable if and only if $n\ge (l-1)(m-1)+1$. 
Further, Voloshin posed the following problem in~\cite{volo2002}.

\begin{prob}[\cite{volo2002}]
	\label{ques1}
For $n,l,m\in\N$ with $n\ge (l-1)(m-1)+1$,
consider a mixed hypergraph $\hyh=(\hyv,\hyc,\hyd)$ of order $n$ with the size of each $\hyc$-edge $l~(\ge 3)$ and the size of each $\hyd$-edge $m~(\ge 2)$. Suppose that $\hyh$ does not contain $\hyk(n,l,m)$ as a subhypergraph. Then what is the minimum number of $\hyc$-edges ($\hyd$-edges) for $\hyh$ to be uncolorable?
\end{prob}

Towards answering Problem~\ref{ques1}, 
in this article we shall focus on a most `non-special' case when $m=l$.  A mixed hypergraph $\hyh=(\hyv,\hyc,\hyd)$ is called a \textit{bi-hypergraph} if $\hyc=\hyd$. Since every edge in a bi-hypergraph is both an edge and a co-edge, 
we can simply write  $\hyh$ as $(\hyv,\hye)$
rather than $(\hyv,\hye,\hye)$, 
and the size of $\hyh$ is equal to
$|\hye|$ rather than $2|\hye|$. 
 A \textit{subhypergraph} of a bi-hypergraph $\hyh=(\hyv,\hye)$ is a bi-hypergraph $\hyh'=(\hyv', \hye')$ with $\hyv'\subseteq \hyv$ and $\hye'\subseteq \hye$, and $\hyh'$ is \textit{induced} if $\hye'=\{e\in \hye: e\subseteq \hyv'\}$.
Obviously, a proper coloring $f$ of a bi-hypergraph $\hyh=(\hyv,\hye)$ satisfies the requirement that $1<|f(e)|<|e|$ for each $e\in\hye$. 
Thus, a bi-hypergraph $\hyh$ is uncolorable if it contains an edge 
with size less than $3$.

Moreover, Tuza and Voloshin~\cite{tuza2000} showed their strong interest in \textit{minimal uncolorable bi-hypergraphs}, i.e., the bi-hypergraphs which are themselves uncolorable but 
all of their proper subhypergraphs are colorable.
It is highlighted in~\cite{tuza2000} that $\hyk(n,r,r)$ is minimal uncolorable if and only if $n= (r-1)^2+1$.
Also, it is not difficult to show that 
	each minimal uncolorable bi-hypergraph $\hyh=(\hyv,\hye)$ 
	is $2$-edge-connected, 
	i.e., both $\hyh$ and $\hyh-e$ for any edge 
	$e$ are connected, where $\hyh-e=(\hyv,\hye\setminus \{e\})$.
The following problem was raised by Tuza and Voloshin~\cite{tuza2000} in 2000 and remains open.

\begin{prob}[\cite{tuza2000}]\label{ques2}
	What is the smallest positive integer $m=m(r)$ for which there exists a minimal uncolorable $r$-uniform bi-hypergraph with $m$ edges?
\end{prob}

\def \muc {{\cal MUC}}
\def \bH {Bi\hyh yp}

For any $n,r\in \N$,  let $\bH(n,r)$
denote the set of $r$-uniform bi-hypergraphs of order $n$,
and let $\muc(n,r)$
denote the set of  minimal uncolorable bi-hypergraphs in $\bH(n,r)$. 
For any $r\ge 3$, if $r\le n\le (r-1)^2$, 
then every bi-hypergraph
in $\bH(n,r)$ is colorable (see \cite{tuza2000}),
implying that $\muc(n,r)=\emptyset$.
If $n\ge (r-1)^2+1$,
we guess that $\muc(n,r)$ is non-empty. 

\begin{con}\label{conj0}
	For any $n,r\in \N$ with $r\ge 3$
	and $n\ge (r-1)^2+1$,
	 the set $\muc(n,r)$ is not empty. 
\end{con}

For any $r\ge 3$ and 
$n\ge  (r-1)^2+1$,
if $\muc(n,r)\ne \emptyset$,  
let $\low{n}{r}$ be the minimum number of edges in a bi-hypergraph 
$\hyh\in \muc(n,r)$;
otherwise, 
 $\low{n}{r}=\infty$. 
 Obviously, $m(r)=\min_{n\in\N} \low{n}{r}$.
 The only known result on $\low{n}{r}$ 
 is that $m(r)\le \low{(r-1)^2+1}{r}= \binom{(r-1)^2+1}{r}$, due to Tuza and Voloshin~\cite{tuza2000}.
 For any $n,r,m\in \N$,  let $\bH(n,r,m)$
 denote the set of bi-hypergraphs in $\bH(n,r)$ of size $m$.
 It is clear that 
 for any $k\in \N$, 
 if every bi-hypergraph 
  in $\bH(n,r,k-1)$ 
  is colorable, then
 $\low{n}{r}\ge k$.

In response to Problems~\ref{ques1} and~\ref{ques2} and Conjecture~\ref{conj0} from a big picture, we will give some sufficient conditions in Section~\ref{sec2}, for the colorability of uniform bi-hypergraphs.
Especially,
we will apply 
probabilistic methods
to show that an $r$-uniform bi-hypergraph is colorable 
when $r$ is relatively large and the edges are sparse or evenly distributed.

\them{thm0}
{
An $r$-uniform bi-hypergraph $\hyh$ is colorable if any of the following requirements holds:
	\begin{enumerate}
	\item $\hyh$ has less than $(r-1)^{r-1}$ edges;
	\item  every edge in $\hyh$ is incident to less than $(r-1)^{r-1}e^{-1}-1$ other edges, 
	 where $e$ ($\approx 2.718281828\ldots$) is the base of natural logarithms.
	\end{enumerate}
}

Note that Theorem~\ref{thm0} determines a lower bound of $m(r)$.

\corr{corrr}
{
For any $r\ge 3$, $m(r)\ge (r-1)^{r-1}$.
}

	Further,
	 Theorem~\ref{thm0} indicates that for uncolorable $r$-uniform bi-hypergraphs $\hyh$ with $r\ge 4$, the size of $\hyh$ should be considerably large, which leads us to focus on the study of $3$-uniform bi-hypergraphs.

In Section~\ref{sec6}, we will study 
Problem~\ref{ques2}  and solve this problem for $r=3$.

\them{thm2}
{$m(3)=10$.
}

In Section~\ref{sec1}, 
we will show that Conjeture~\ref{conj0} holds when $r=3$
and find an upper bound for 
$\low{n}{3}$, as stated below.

\them{thm1}
{
	For any $n\in\N$ with $n\ge 6$, $\muc(n,3)\ne \emptyset$, and
	$\low{n}{3}\le \frac{7n}{3}-4$.
}

We will prove Theorem~\ref{thm1} 
by constructing a  minimal
uncolorable $3$-uniform
bi-hypergraph of order $n$ for every $n\ge 6$, each of which does not contain $\hyk(5,3,3)$ as subhypergraphs.

For a bi-hypergraph 
$\hyh=(\hyv,\hye)$ and $v\in \hyv$,
let $\hye_{\hyh}(v)$ denote the set of edges $e\in \hye$ with $v\in e$, 
let $N_{\hyh}(v)$ (or simply 
$N(v)$) denote the set of vertives 
$u\in \hyv\setminus \{v\}$
such that $u\in e$ for some $e\in 
\hye_{\hyh}(v)$, 
let $N_{\hyh}[v]=N_{\hyh}(v)\cup\{v\}$,
and let $d_{\hyh}(v)$ (or simply 
$d(v)$), called 
the \textit{degree} of $v$ in $\hyh$, 
denote the size of $\hye_{\hyh}(v)$.
A set $\hyv_0\subseteq \hyv$ is \textit{independent} if $e\not\subseteq \hyv_0$ for all $e\in \hye$.

\section{Colorable uniform bi-hypergraphs
	\label{sec2}
}

In this section, we bring up with some sufficient conditions for proving bi-hypergraphs being colorable.

	We first give the proof of Theorem~\ref{thm0}, which focuses more on the perspective of edge distribution. The following Lemma~\ref{LLL}, one of the most powerful tools from probabilistic methods, will be applied.


\begin{lem}[Lov\'{a}sz local lemma,~\cite{Erdos1975}]\label{LLL}
	Let $A_1,A_2,\dots,A_m$ be events in an arbitrary probability space. Suppose that the probability of each of the $m$ events is at most $p$, and suppose that each event $A_i$ is mutually independent of all but at most $b$ of the other events $A_j$. If $ep(b + 1) <1$, where $e$ is the base of natural logarithms, then with positive probability none of the events $A_i$ holds.
\end{lem}

\vspace{0.2 cm}

\noindent {\it Proof of Theorem~\ref{thm0}}: \
	Assume that $\hyh=(\hyv,\hye)$, where $|\hyv|=n$ and $\hye=\{e_1,e_2,\dots,e_m\}$. It suffices to show that there is a mapping $c$: $\hyv \rightarrow [r-1]$ 
		such that $|c(e_i)|\ge 2$ for all $i\in[m]$.
		
		Let $\Omega$ be the set of 
		mappings $f$ from $\hyv$ to $[r-1]$ by choosing $f(v)$ for each $v\in \hyv$ randomly and independently. Then $\Omega$ is a uniformly distributed probability space. 
		
		For each $i\in[m]$, let $A_i$ be the event that edge $e_i$ is not properly colored, i.e., $A_i$ consists of the mappings $c$ in $\Omega$ such that $|c(e_i)|=1$. Then 
		\equ{eq2-1}
		{
		P(A_i)=\frac{(r-1)^{n-r+1}}{(r-1)^n}=(r-1)^{-r+1}.
	}
For (i), where $m<(r-1)^{r-1}$, we have
		\equ{eq2-2}
	{
	P(\bigcup_{i=1}^m A_i)\le \sum_{i=1}^m P(A_i) <(r-1)^{r-1} \cdot (r-1)^{-r+1}=1,
}
	which implies that there is a proper coloring of $\hyh$ among set $\Omega$.
	
For (ii), since every event $A_i$ is mutually independent of all but at most $\lceil (r-1)^{r-1}e^{-1}-2\rceil$ of the other events $A_j$, the existence of a proper coloring of $\hyh$ with at most $r-1$ colors is guaranteed by Lemma~\ref{LLL}.

	Hence the result is proven.
	\proofend

Then the following corollary is immediate.
	
	\corr{coro1}
	{
For any $r$-uniform  bi-hypergraph
$\hyh$,  if $d(v)\le \frac{(r-1)^{r-1}e^{-1}-2}{r}$ holds 
for each vertex $v$ in $\hyh$, then $\hyh$ is colorable.
}

In the following, we propose some sufficient conditions regarding vertex partitions that will play an important role in the study of colorability of bi-hypergraphs in
the next section.

	\lemm{le5-1}
	{
		For any $r\ge 3$ and $r$-uniform 
		bi-hypergraph $\hyh=(\hyv, \hye)$, 
		if $\hyv$ has a partition 
		$\hyv_1,\dots,\hyv_k$ such that 
		each $\hyv_i$ is an independent 
		set of $\hyh$
		and for each edge $e\in \hye$,
		$e\subseteq \hyv_{i_1}\cup \cdots\cup \hyv_{i_{q}}$
		holds for some 
		numbers $i_1,\dots, i_{q}$ in $[k]$,
		where $q<r$, 
		then $\hyh$ is colorable. 
	}
	
	\proof Let $c$ be the mapping from $\hyv$ to $\N$ defined by 
	$c(u)=i$ for all $u\in \hyv_i$
	and all $i\in [k]$. 
	By the given conditions, 
	 $2\le |c(e)|\le r-1$ holds for each $e\in \hye$.
	Thus, 
	$c$ is a proper coloring of $\hyh$.
	\proofend 

Two special cases of Lemma~\ref{le5-1}
are given below.

\corr{cor5-1}
{
For any $r\ge 3$ and $r$-uniform 
bi-hypergraph $\hyh=(\hyv, \hye)$, 
if $\hyv$ has a partition 
$\hyv_1,\dots,\hyv_k$, where $k\le r-1$ and  each $\hyv_i$ is independent,
then $\hyh$ is colorable. 
}

\corr{cor5-2}
{
	For any $r\ge 3$ and $r$-uniform 
	bi-hypergraph $\hyh=(\hyv, \hye)$, 
	if $\hyv$ has a partition 
	$\hyv_1,\dots,\hyv_r$
	such that $\hyv_i$ is independent
	for all $i\in [r]$, 
	$\hyv_r=\{w\}$
	and for each $e\in \hye$ with $w\in e$, 
	$e\cap \hyv_i=\emptyset$ holds for some 
	$i\in [r-1]$, 
	then $\hyh$ is colorable. 
}

	Applying Corollary~\ref{cor5-1}, we can prove the following result.

	\lemm{prop5}
	{
		For any $r\ge 3$ and 
		$\hyh=(\hyv, \hye)\in \bH(2r,r)$, 
		if 
		$|\hye|<{2r\choose r}/2$,
		then $\hyh$ is colorable. 
	}
	
	\proof
	As $|\hyv|=2r$, $\hyv$ has exactly
	${2r\choose r}$ $r$-element subsets.
	Thus, ${\hyv\choose r}$ can be partitioned into 
	${2r\choose r}/2={2r-1\choose r-1}$ 
	subsets:
	$
	S_1, S_2,\dots, S_s,
	$
	where $s={2r\choose r}/2$ and 
	$S_i=\{e_i,e_i'\}$ with 
	$e_i\cup e'_i=\hyv$
	for each $i\in [s]$. 
	
	Since $|\hye|<{2r\choose r}/2$, 
	there exists a set $S_i=\{e_i,e_i'\}$
	such that both $e_i$ and $e'_i$ 
	are not edges in $\hyh$. 
	Since $e_i$ and $e_i'$ form a partition
	of $\hyv$, 
	the conclusion follows from 
	Corollary~\ref{cor5-1}.
	\proofend

In fact, Lemma~\ref{prop5} determines a lower bound on $\low{6}{3}$.
\corr{prop6-3}
{
Any $\hyh=(\hyv,\hye) \in\bH(6,3)$ with $|\hye|\le 9$ is colorable. Hence $\low{6}{3}\ge 10$.
}

The next lemma is quite useful
as it makes it possible to prove colorability by identifying vertices.

\lemm{lem15}
{
	Let $\hyh$ be a bi-hypergraph with two non-adjacent vertices $u$ and $v$, and let $\hyh\cdot uv$ be the bi-hypergraph obtained from $\hyh$ by identifying vertices $u$ and $v$. 
	If $\hyh\cdot uv$ is colorable, then  $\hyh$ is also colorable.
}

\proof
Assume that $\hyh=(\hyv,\hye)$ and $\hyh \cdot uv=(\hyv',\hye')$, where $\hyv'\setminus \hyv=\{w\}$.

Let $c$ be a proper coloring of $\hyh \cdot uv$. Then let $c'$ be the mapping from $\hyv$ to $\N$ such that $c'(x)=c(x)$ for any $x\in \hyv\cap \hyv'$ and $c'(u)=c'(v)=c(w)$. We shall show that $c'$ is a proper coloring of $\hyh$.

Since $u$ and $v$ are not adjacent in $\hyh$, for each edge $e$ in $\hye$, $|e\cap \{u,v\}|\le 1$. The edges $e$ with $e\cap \{u,v\}=\emptyset$ are obviously properly colored. For the remaining edges $e$ such that $|e\cap \{u,v\}|=1$, let $e'$ be the corresponding edge of $e$ in $\hyh \cdot uv$, i.e., $e'=(e\setminus\{u,v\})\cup \{w\}$. As a result, $|e'|=|e|$ and $|\{c(e')\}|=|\{c'(e)\}|$, which implies that $e$ is properly colored. Hence $\hyh$ is colorable.
\proofend

\section{Proof of Theorem~\ref{thm2}
\label{sec6}}

In this section,
we will prove Theorem~\ref{thm2}. Due to the fact that $\hyk(5,3,3)$ is minimal uncolorable and Corollary~\ref{prop6-3},
it remains to show that every $3$-uniform bi-hypergraph with order at least seven and size up to nine is colorable.

Recall that $\bH(n,r,m)$ is the set of  $r$-uniform bi-hypergraphs 
	with $n$ vertices and $m$ edges.
By the Handshaking Lemma, the following observation is direct.

\obsv{obsv6}
{
	For any $\hyh\in \bH(n,r,m)$,
	 there is a vertex $v$ in $\hyh$
	with $d(v)\le  \lfloor \frac{rm}{n}\rfloor$. 
}

As a special case of Observation~\ref{obsv6}, if $n\ge 7$ and $m\le 9$, then every bi-hypergraph 
in $\bH(n,3,m)$ has a vertex 
$v$ with $d(v)\le 3$.

\lemm{lem16}
{
	For any $n,m,r\in \N$ with $n\ge \lfloor \frac{rm}{n+1}\rfloor (r-1)+1$,
	if all bi-hypergraphs in $\bH(n,r,m)$ are colorable, 
	then all bi-hypergraphs in $\bH(n+1,r,m)$ are colorable. 
}

\proof
Let $\hyh=(\hyv,\hye)\in \bH(n+1,r,m)$, where $n\ge  \lfloor \frac{rm}{n+1}\rfloor (r-1) + 1$.
By Observation~\ref{obsv6}, 
$d(u)\le \lfloor \frac{rm}{n+1}\rfloor$
holds for some  $u\in \hyv$. 
Since $\hyh$ is $r$-uniform, 
	\equ{eq3-1}
{
|N_{\hyh}[u]|\le \left \lfloor \frac{rm}{n+1}\right \rfloor \cdot (r-1)+1\le n,
}
which implies that
there exists $v\in \hyv$ such that 
$u$ and $v$ are not adjacent in $\hyh$. 
Thus, $\hyh\cdot uv\in \bH(n,r,m)$. 
By Lemma~\ref{lem15}, 
whenever $\hyh \cdot uv$ is colorable, $\hyh$ is also colorable. 
Hence the result holds.
\proofend

The special case of Lemma~\ref{lem16} for $r=3$ is given below. 

\corr{co-lem16}
{
	For any $n,m\in \N$ with $n\ge 7$ and $m\le 9$, 
	if all bi-hypergraphs in $\bH(n,3,m)$ are colorable, 
	then all bi-hypergraphs in $\bH(n+1,3,m)$ are colorable. 
}

Corollary~\ref{co-lem16} indicates that, 
 in order to prove Theorem~\ref{thm2}, 
 we need only to prove that every bi-hypergraph in $\bH(7,3)$ with size at most nine is colorable.

\begin{figure}[!ht]
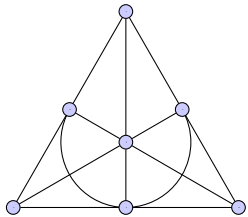

	\tikzp{1.5}
	{
		
		\draw (0,0) +(150:0.577) arc (150:360:0.577) ;
		\draw (0,0) +(0:0.577) arc (0:30:0.577) ;
		
		\foreach \place/\y in {{(0,1.1546)/1},{(0.5,0.2886)/2},{(-0.5,0.2886)/3},
			{(0,0)/4}, {(1,-0.5773)/5},{(0,-0.5773)/6},{(-1,-0.5773)/7}}
		\node[cblue] (b\y) at \place {};
		
		\filldraw[black] (b1) circle (0pt)node[anchor=north] {};
		\filldraw[black] (b2) circle (0pt)node[anchor=north] {};
		\filldraw[black] (b3) circle (0pt)node[anchor=north] {};
		\filldraw[black] (b4) circle (0pt)node[anchor=south] {};
		\filldraw[black] (b5) circle (0pt)node[anchor=south] {};
		\filldraw[black] (b6) circle (0pt)node[anchor=south] {};
		\filldraw[black] (b7) circle (0pt)node[anchor=south] {};		
		
		\draw (b7) -- (b6) -- (b5) -- (b2) -- (b1);
		\draw (b2) -- (b4) -- (b7) -- (b3) -- (b1) -- (b4) -- (b6);		
		\draw (b3) -- (b4) -- (b5);
		
	}
	
	\caption{Bi-hypergraph $\hyf\hyp$, where every line segment represents an edge}
	
	\label{fig1-3}
\end{figure}

A key structure that will emerge is a special $3$-uniform bi-hypergraph $\hyf\hyp$,
which is also known as the {\it Fano plane hypergraph}.
As shown in Figure~\ref{fig1-3},
$\hyf\hyp$ can be represented by 
a diagram that usually represents
the Fano matroid $F_7$ in which the $3$-element circuits correspond to edges in 
$\hyf\hyp$.

\lemm{lem17}
{
Any $\hyh=(\hyv,\hye) \in\bH(7,3)$ with $|\hye|\le 9$ is colorable. Hence $\low{7}{3}\ge 10$.
}

\proof 
Suppose there exists 
$\hyh=(\hyv, \hye)\in \bH(7,3)$
with $|\hye|\le 9$ such that 
$\hyh$ is uncolorable.
Now we are going to prove the following 
claims.

\Clm{cl1}{For any $x\in \hyv$, $N(x)=\hyv\setminus\{x\}$.}
\proof
Suppose that $x,y$ are not adjacent in $\hyh$. As a result, $\hyh\cdot xy$ is in $\bH(6,3)$ with size at most $9$. By Corollary~\ref{prop6-3}, $\hyh\cdot xy$ is colorable, which implies that $\hyh$ is colorable by Lemma~\ref{lem15}.
\claimend


\Clm{cl2}{$d(u)=3$ for some $u\in \hyv$.
}

\proof
By Observation~\ref{obsv6}, 
there exists $u\in \hyv$ with $d(u)\le 3$.
Since $\hyh$ is $3$-uniform and $|\hyv\setminus\{u\}|=6$, 
by Claim~\ref{cl1}, $d(u)\ge 3$.
Thus, $d(u)=3$
and Claim~\ref{cl2} follows.
\claimend

By Claim~\ref{cl2}, the subhypergraph of $\hyh$ induced by $\hyv\setminus \{u\}$, denoted by $\hyh-u$, has order 6 and size at most 6.
Similar to the proof of
Lemma~\ref{prop5},  
$\hyv\setminus \{u\}$ has a partition 
 $V_1=\{v_i: 1\le i\le 3\}$ 
 and $V_2=\{v_i: 4\le i\le 6\}$ such that $V_i\notin \hye$ for both $i=1,2$.

\Clm{cl3}{For each $i=1,2$, there is  one edge $e_i\in \hye$  such that
	$u\in e_i\subseteq V_i\cup \{u\}$.}

\proof
Suppose $\hyh$ has no edge $e$ 
such that 
$u\in e\subseteq V_1\cup \{u\}$. 
Then $V_1\cup \{u\}$ 
is independent in $\hyh$,
implying that $\hyv$ can be partitioned
into two independent sets
$V_1\cup \{u\}$ and $V_2$. 
By Corollary~\ref{cor5-1}, 
$\hyh$ is colorable, 
a contradiction to the assumption.
\claimend

\begin{figure}[!ht]
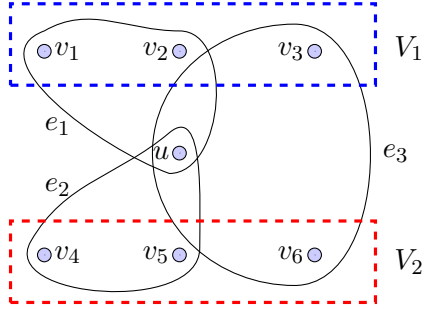

	\tikzp{0.9}
	{
		\foreach \place/\y in {{(-2,2.5)/1},{(0,1)/2},{(2,2.5)/3},
			{(-2,-0.5)/4}, {(0,-0.5)/5},{(2,-0.5)/6},{(0,2.5)/7}}
		\node[cblue] (b\y) at \place {};
		
		\filldraw[black] (b1) circle (0pt)node[anchor=west] {$v_1$};
		\filldraw[black] (b2) circle (0pt)node[anchor=east] {$u$};
		\filldraw[black] (b3) circle (0pt)node[anchor=east] {$v_3$};
		\filldraw[black] (b4) circle (0pt)node[anchor=west] {$v_4$};
		\filldraw[black] (b5) circle (0pt)node[anchor=east] {$v_5$};
		\filldraw[black] (b6) circle (0pt)node[anchor=east] {$v_6$};
		\filldraw[black] (b7) circle (0pt)node[anchor=east] {$v_2$};

		\begin{scope}[fill opacity=0.2]

			\draw ($(b2)+(-0.15,-0.3)$) 
			to[out=160,in=200] ($(b1) + (0,0.3)$) 
			to[out=30,in=180] ($(b7) + (0,0.3)$)
			to[out=0,in=00] ($(b2) + (-0.15,-0.3)$);    
			
			\draw ($(b2)+(-0.15,0.3)$) 
			to[out=220,in=60] ($(b4) + (-0.2,0.2)$) 
			to[out=240,in=270] ($(b5) + (0.3,0.1)$)
			to[out=90,in=40] ($(b2) + (-0.15,0.3)$);       
			
			\draw ($(b2)+(-0.4,0)$) 
			to[out=90,in=180] ($(b3) + (-0.3,0.4)$) 
			to[out=0,in=0] ($(b6) + (-0.3,-0.45)$)
			to[out=180,in=270] ($(b2) + (-0.4,0)$);                                                                        
		\end{scope}
		
		\draw[blue, dashed, very thick] (-2.5, 3.2) rectangle (2.9,2);
		\draw[red, dashed, very thick] (-2.5, -1.2) rectangle (2.9,0);
		
		\node [style=none] (cap1) at (3.4, 2.55) {$V_1$};			
		\node [style=none] (cap2) at (3.4, -0.6) {$V_2$};	
		
	\node [style=none] (cap1) at (-1.8, 1.4) {$e_1$};	
	
	\node [style=none] (cap1) at (-1.8, 0.5) {$e_2$};	
		
	\node [style=none] (cap1) at (3.2, 1.0) {$e_3$};	
				
	}
	\caption{Edges containing vertex $u$}
	\label{fig6}
\end{figure}

Assume that $e_1=\{u,v_1,v_2\}$
and $e_2=\{u,v_4,v_5\}$.
By Claim~\ref{cl2}, $d(u)=3$,
and by Claim~\ref{cl1},
$N(u)=\hyv\setminus \{u\}=V_1\cup V_2$.
Thus,  another edge in $\hyh$ containing $u$ is $e_3=\{u,v_3,v_6\}$,
 as shown in Figure~\ref{fig6}.

\Clm{cl4}{For any $v$ in $V_i$, where $i\in[2]$, there exists $e\in \hye$
	such that $v\in e\subseteq V_{3-i}\cup \{v\}$. 
 }

\proof
We need only to consider the cases
$v=v_1$ and $v=v_3$.

If $\hyh$ has no edge $e$ 
such that $v_1\in e\subseteq V_{2}\cup \{v_1\}$, then $\hyv$ can be partitioned 
into two independent sets
$V_2\cup \{v_1\}$ and $\{u,v_2,v_3\}$,
implying that 
$\hyh$ is colorable
by Corollary~\ref{cor5-1}.

If $\hyh$ has no edge $e$ 
such that $v_3\in e\subseteq V_{2}\cup \{v_3\}$, then $\hyv$ can be partitioned 
into three independent sets
$\{v_j: 3\le j\le 6\}$, $\{v_1, v_2\}$
and $\{u\}$.
Since $e_1,e_2$ and $e_3$ are the only edges in $\hyh$ 
containing $u$
and 
either $e_i\cap \{v_j: 3\le j\le 6\}=\emptyset$ 
or $e_i\cap \{v_1,v_2\}=\emptyset$ 
for each $i\in [3]$, 
$\hyh$ is colorable
by Corollary~\ref{cor5-2}.

Both cases lead to a contradiction, 
and thus Claim~\ref{cl4} holds.
\claimend

Since $|\hye|\le 9$ and $d(u)=3$, 
$\hyh-u$ contains at most six edges.
By Claim~\ref{cl4},  
$\hyh-u$ has at least six edges.
Thus, the following claim follows. 

\Clm{cl4-0}{$|\hye|=9$. }

By Claims~\ref{cl2} and~\ref{cl4-0}, 
$\hyh-u$ has exactly six edges.
Furthermore, we have the following stronger version of Claim~\ref{cl4}.

\Clm{cl4'}{For any $v$ in $V_i$, where $i\in[2]$,  there exists a unique edge 
	$e\in \hye$ such that 
	$v\in e\subseteq V_{3-i}\cup \{v\}$.
}

Note that  $e_1=\{u,v_1,v_2\}$, $e_2=\{u,v_4,v_5\}$
and $e_3=\{u,v_3,v_6\}$. 

\Clm{cl5}{For any 
	$I\subseteq [6]$ with $|I|=3$,
	if $|\{v_j: j\in I\}\cap e_i|=1$ holds 
	for all $i\in [3]$, 
	then, either $\{v_j:j\in I\}\in\hye$ or $\{v_j:j\in I'\}\in\hye$,
	where $I'=[6]\setminus I$.
}

\proof For any 
$I\subseteq [6]$ with $|I|=3$,
if the conclusion fails, then 
$\hyv$ can be partitioned 
into two independent sets in $\hyh$: 
$\{u\}\cup \{v_j: j\in I\}$
and $\{v_j: j\in I'\}$,
implying that $\hyh$ is colorable
by Corollary~\ref{cor5-1}, 
a contradiction.
\claimend

There are exacly four 
non-ordered pairs $(I,I')$, 
where $I'=[6]\setminus I$, 
such that $I\subseteq [6]$, 
$|I|=3$
and $|\{v_j: j\in I\}\cap e_i|=1$ holds 
for all $i\in [3]$:
$$
(\{1,3,4\},\{2,5,6\}); (\{1,3,5\},\{2,4,6\}); (\{1,4,6\},\{2,3,5\}); (\{1,5,6\},\{2,3,4\}).
$$ 

\Clm{cl6}{$\hyh$ has a spanning subhypergraph $\hyh_0=(\hyv,\hye_0)$,
	as shown in Figure~\ref{fig1-4} 
	(a) or (b),
 which is isomorphic to $\hyf\hyp$.
}
\proof
First consider the case 
$\{v_1,v_3,v_4\}\in\hye$. 
Then by Claim~\ref{cl4'}, $\{v_2,v_3,v_4\}\notin \hye$, which implies that $\{v_1,v_5,v_6\}\in\hye$. Similarly, $\{v_1,v_4,v_6\}\notin \hye$, and thus $\{v_2,v_3,v_5\}\in\hye$. Finally, $\{v_1,v_3,v_5\}\notin \hye$ and so $\{v_2,v_4,v_6\}\in\hye$. 
Now it is easy to check $\hyh_0=(\hyv,\hye_0)$ is a subhypergraph of $\hyh$ that is isomorphic to $\hyf\hyp$, where 
$$
\hye_0=\{e_1,e_2,e_3\}
\cup 
\{ 
\{v_1,v_3,v_4\}, \{v_1,v_5,v_6\},  \{v_2,v_3,v_5\}, \{v_2,v_4,v_6\}\},
$$
as shown in Figure~\ref{fig1-4} (a).

If $\{v_2,v_5,v_6\}\in\hye$, 
then it can be deduced similarly that 
$\hyh$ has a spanning subhypergraph 
$\hyh_0=(\hyv,\hye_0)$  
isomorphic to $\hyf\hyp$, where 
$$
\hye_0=\{e_1,e_2,e_3\}
\cup 
\{ 
\{v_2,v_5,v_6\}, \{v_2,v_3, v_4\},  \{v_1,v_4, v_6\}, \{v_1,v_3, v_5\}\},
$$
as shown in Figure~\ref{fig1-4} (b).
\claimend

\begin{figure}[!ht]
	\tikzp{1.5}
	{
		\draw (0,0) +(150:0.577) arc (150:360:0.577) ;
		\draw (0,0) +(0:0.577) arc (0:30:0.577) ;
		
		\foreach \place/\y in {{(0,1.1546)/1},{(0.5,0.2886)/2},{(-0.5,0.2886)/3},
			{(0,0)/4}, {(1,-0.5773)/5},{(0,-0.5773)/6},{(-1,-0.5773)/7}}
		\node[cblue] (b\y) at \place {};
		
		\filldraw[black] (b1) circle (0pt)node[anchor=south] {$v_6$};
		\filldraw[black] (b2) circle (0pt)node[anchor=west] {$v_5$};
		\filldraw[black] (b3) circle (0pt)node[anchor=east] {$v_2$};
		\filldraw[black] (b4) circle (0pt)node[anchor=south west] {$u$};
		\filldraw[black] (b5) circle (0pt)node[anchor=north] {$v_1$};
		\filldraw[black] (b6) circle (0pt)node[anchor=north] {$v_3$};
		\filldraw[black] (b7) circle (0pt)node[anchor=north] {$v_4$};		
		
		\draw (b7) -- (b6) -- (b5) -- (b2) -- (b1);
		\draw (b2) -- (b4) -- (b7) -- (b3) -- (b1) -- (b4) -- (b6);		
		\draw (b3) -- (b4) -- (b5);
		\draw (4,0) +(150:0.577) arc (150:360:0.577) ;
		\draw (4,0) +(0:0.577) arc (0:30:0.577) ;
		
		\foreach \place/\y in {{(4,1.1546)/1},{(4.5,0.2886)/2},{(3.5,0.2886)/3},
			{(4,0)/4}, {(5,-0.5773)/5},{(4,-0.5773)/6},{(3,-0.5773)/7}}
		\node[cblue] (c\y) at \place {};
		
		\filldraw[black] (c1) circle (0pt)node[anchor=south] {$v_3$};
		\filldraw[black] (c2) circle (0pt)node[anchor=west] {$v_4$};
		\filldraw[black] (c3) circle (0pt)node[anchor=east] {$v_1$};
		\filldraw[black] (c4) circle (0pt)node[anchor=south west] {$u$};
		\filldraw[black] (c5) circle (0pt)node[anchor=north] {$v_2$};
		\filldraw[black] (c6) circle (0pt)node[anchor=north] {$v_6$};
		\filldraw[black] (c7) circle (0pt)node[anchor=north] {$v_5$};		
		
		\draw (c7) -- (c6) -- (c5) -- (c2) -- (c1);
		\draw (c2) -- (c4) -- (c7) -- (c3) -- (c1) -- (c4) -- (c6);		
		\draw (c3) -- (c4) -- (c5);
	}
	
	\centerline{(a) \hspace{5 cm} (b)}

	\caption{A bi-hypergraph $\hyh_0$
		isomorphic to  $\hyf\hyp$}
	
	\label{fig1-4}
\end{figure}

Applying the above claims,  
we are now able to prove 
that $\hyh$ is colorable. 

\Clm{cl7}{$\hyh$ is colorable.}

\proof 
By  Claim~\ref{cl4'}, 
$\hyh$ has edges $e$ and $e'$ such that  
\equ{eq6-1}
{
	v_3\in e\subseteq V_2\cup \{v_3\}
	\quad \mbox{and}\quad 
	v_6\in e'\subseteq V_1\cup \{v_6\}.
}
Clearly, both $e$ and $e'$ are not edges in $\hyh_0$.
Since $|\hye|=9$ by Claim~\ref{cl4-0}, 
we have $\hye=\hye_0\cup \{e,e'\}$.

Since $\hye=\hye_0\cup \{e,e'\}$,
the set $\{v_1,v_2,v_4,v_5\}$ is independent in $\hyh$.
It follows that 
$\hyv$ can be partitioned into 
three independent sets 
$\{v_1,v_2,v_4,v_5\}$, $\{v_3,v_6\}$
and $\{u\}$. 
Clearly, for each $i\in [3]$,
either $e_i\cap \{v_1,v_2,v_4,v_5\}=\emptyset$
or $e_i\cap \{v_3,v_6\}=\emptyset$.
Since $e_1,e_2$ and $e_3$ are the only edges 
in $\hyh$ containing $u$, 
by Corollary~\ref{cor5-2}, 
$\hyh$ is colorable. 
\claimend

Claim~\ref{cl7}  contradicts with the assumption. Thus, the result is proven.
\proofend

Theorem~\ref{thm2} follows directly 
from Lemma~\ref{lem17}, 
Corollaries~\ref{prop6-3}
and~\ref{co-lem16} 
and the fact that $\hyk(5,3,3)$ is minimal uncolorable.

\section{Proof of Theorem~\ref{thm1}
\label{sec1}}

In this section, we shall construct a number of minimal uncolorable $3$-uniform bi-hypergraphs that 
do not contain $\hyk(5,3,3)$.
The key step in our construction is a repetitive structure 
given in Subsection~4.1 which possesses nice properties.

\subsection{A $3$-uniform bi-hypergraph
	of order $3k$ 
and size $7k-6$}

\defn{def1}
{Let $k\ge 2$, and let 
$V_1,V_2,\dots,V_k$ be pairwise disjoint  sets, where 
$V_i= \{\vtx{i}{j}:j\in [3]\}$
for all $i\in [k]$.
We define the bi-hypergraph $\hyh_{k}$ as follows:
\begin{enumerate}
\item 
for all $q\in[k-1]$, $V_q\times V_{q+1}$ is 
the bi-hypergraph $(\hyv, \hye)$, where
$\hyv=V_q\cup V_{q+1}$
and
$$
\hye=
\{ \{\vtx{q+t}{j}:j\in [3]\}
:t=0,1 \}
\cup 
\left \{ 
\{ \vtx{q+1}{j},\vtx{q}{j},\vtx{q}{j+t}\}: j=1,2,3,t=1,2
\right \},
$$
where $\vtx{q}{4}=\vtx{q}{1}$ and $\vtx{q}{5}=\vtx{q}{2}$;
\item $\hyh_2$ is the bi-hypergraph $V_1\times V_2$, and for all $k\ge 3$,
$\hyh_{k}$ is the bi-hyergraph 
obtained from  
$\hyh_{k-1}$ 
by adding all the vertices in $V_{k}$ and all the edges in $V_{k-1}\times V_k$.
\end{enumerate} 
}

The structure of $\hyh_2$ can be seen from its vertex-edge incidence bipartite graph shown in Figure~\ref{fig1-1},
where 
the  {\it vertex-edge incidence graph} (or simply {\it the incidence graph})
of a bi-hypergraph 
$\hyh=(\hyv, \hye)$ is a bipartite graph
$G$  
with a bipartition $(\hyv, \hye)$ 
such that $v\in \hyv$ is adjacent to 
$e\in \hye$ in $G$ if and only if 
$v\in e$. 

\begin{figure}[!ht]
	\tikzp{1.2}
	{
		\foreach \place/\y in {{(-2.5,-1.2)/1},{(-0.5,-1.2)/2},{(1.5,-1.2)/3},
			{(-2.5,1.2)/4}, {(-0.5,1.2)/5},{(1.5,1.2)/6}}
		\node[cblue] (b\y) at \place {};
		
		\filldraw[black] (b1) circle (0pt)node[anchor=north] {$\vtx{1}{1}$};
		\filldraw[black] (b2) circle (0pt)node[anchor=north] {$\vtx{1}{2}$};
		\filldraw[black] (b3) circle (0pt)node[anchor=north] {$\vtx{1}{3}$};
		\filldraw[black] (b4) circle (0pt)node[anchor=south] {$\vtx{2}{1}$};
		\filldraw[black] (b5) circle (0pt)node[anchor=south] {$\vtx{2}{2}$};
		\filldraw[black] (b6) circle (0pt)node[anchor=south] {$\vtx{2}{3}$};

		\foreach \place/\y in {{(-4,0)/1},{(-3,0)/2},{(-2,0)/3},{(-1,0)/4},{(0,0)/5},{(1,0)/6},{(2,0)/7},{(3,0)/8}}
		\node[cblue] (a\y) at \place {};

		

		\filldraw[black] (a1) circle (2pt)node[anchor=east] {$e_1$};
		
		\filldraw[black] (a2) circle (2pt)node[anchor=east] {$e_2$};
		
		\filldraw[black] (a3) circle (2pt)node[anchor=east] {$e_3$};
		
		\filldraw[black] (a4) circle (2pt)node[anchor=east] {$e_4$};
		
		\filldraw[black] (a5) circle (2pt)node[anchor=west] {$e_5$};
		
		\filldraw[black] (a6) circle (2pt)node[anchor=west] {$e_6$};
		
		\filldraw[black] (a7) circle (2pt)node[anchor=west] {$e_7$};
		
		\filldraw[black] (a8) circle (2pt)node[anchor=west] {$e_8$};
		
		
		
		\path[-] (a1) edge (b1) edge (b2) edge (b3);
		
		
		
		\path[-] (a2) edge (b4) edge (b1) edge (b2);

		\draw (b4) -- (a3) -- (b1);
		\draw (a3) -- (b3);
		
		\draw (b5) -- (a4) -- (b2);
		\draw (a4) -- (b3);
		
		\draw (b5) -- (a5) -- (b2);
		\draw (a5) -- (b1);
		
		\draw (b6) -- (a6) -- (b3);
		\draw (a6) -- (b1);
		
		\draw (b6) -- (a7) -- (b3);
		\draw (a7) -- (b2);
		
		\draw (b4) -- (a8) -- (b5);
		\draw (a8) -- (b6);
		
		
		
	}
	
	\caption{The incident graph of 
		$\hyh_2$}
	
	\label{fig1-1}
\end{figure}

	Note that $\hyh_k$ has exactly $3k$ vertices 
	and $7k-6$ edges.
	For any $2\le s\le k-1$,
	the subhypergraph of $\hyh_k$ 
	induced by $V_{s-1}\cup V_s$ 
	is isomorphic to the one 
	induced by $V_{s}\cup V_{s+1}$,
	with an isomorphism 
	$f: V_{s-1}\cup V_s
	\rightarrow V_{s}\cup V_{s+1}$
	defined by 
	$f(v_{j,t})=v_{j+1,t}$ 
	for all $j\in \{s-1,s\}$ 
	and $t\in [3]$.

Two basic properties on proper colorings of $\hyh_k$ are given below.

	\obsv{obsv1}
	{
		\begin{enumerate}
			\item For any proper coloring $c$ of $\hyh_2$, 
			there exists $j_1, j_2\in \N$ with 
			$1\le j_1<j_2\le 3$ such that 
			$$
			c(\vtx{1}{j_1})=c(\vtx{1}{j_2})=
			c(\vtx{2}{z})=a, \quad
			c(\vtx{2}{j_1})=c(\vtx{2}{j_2})=
			c(\vtx{1}{z})=b
			$$
			for some $a,b\in \N$, where $a\ne b$
			and $z\in [3]\setminus \{j_1,j_2\}$;
			
			\item for any $j_1, j_2\in \N$ with 
			$1\le j_1<j_2\le 3$,  
			$\hyh_2$ has a proper coloring $c$ such that 
			$c(\vtx{1}{j_1})=c(\vtx{1}{j_2})$. 
			
		\end{enumerate}  
	}
\proof
Let $c$ be a proper coloring of
$\hyh_2$. 
Since $\{v_{1,i}: i\in [3]\}$
is an edge in $\hyh_2$, 
there exists $j_1,j_2\in\N$ with $1\le j_1<j_2\le 3$ such that 
$c(\vtx{1}{j_1})=c(\vtx{1}{j_2})=a$ and $c(\vtx{1}{z})=b$
	for some $a,b\in \N$, where $a\ne b$
			and $z\in [3]\setminus \{j_1,j_2\}$.

Due to edges $\{\vtx{2}{j_1}, \vtx{1}{j_1}, \vtx{1}{j_2}\}$ and $\{\vtx{2}{j_1}, \vtx{1}{j_1}, \vtx{1}{z}\}$, we have $c(\vtx{2}{j_1})=b$. 
Similarly, $c(\vtx{2}{j_2})=b$ due to edges $\{\vtx{2}{j_2}, \vtx{1}{j_2}, \vtx{1}{j_1}\}$ and $\{\vtx{2}{j_2}, \vtx{1}{j_2}, \vtx{1}{z}\}$. 
Then, due to  edges $\{\vtx{2}{z}, \vtx{1}{z},\vtx{1}{j_1}\}$ and $\{\vtx{2}{j_1},\vtx{2}{j_2},\vtx{2}{z}\}$, we have $c(\vtx{2}{z})=a$. Hence (i) holds.

Then (ii) follows from the coloring given in (i) for any $j_1, j_2\in \N$ with 
			$1\le j_1<j_2\le 3$.
\proofend

\obsv{obsv2}
{
Let  $k\ge 3$. 
For any proper coloring $c$ of
$\hyh_k$, 
and any  $i_1, i_2\in [k]$ 
with $i_1-i_2\equiv 0 \pmod{2}$, 
$c(\vtx{i_1}{j})=c(\vtx{i_2}{j})$ holds
for all $j\in [3]$. 
}

\proof
We only need to show that $c(\vtx{1}{j})=c(\vtx{3}{j})$ holds for all $j\in [3]$.

Without loss of generality, we assume that $c(\vtx{1}{1})=c(\vtx{1}{2})=a$ and $c(\vtx{1}{3})=b$ for some $a,b\in \N$, where $a\ne b$.
Then, according to Observation~\ref{obsv1} (i),
$c(\vtx{2}{1})=c(\vtx{2}{2})=b$ and $c(\vtx{2}{3})=a$ hold.
Applying Observation~\ref{obsv1} (i) again, we have that 
$c(\vtx{3}{1})=c(\vtx{3}{2})=a$ and $c(\vtx{3}{3})=b$.

Hence the result is proven.
\proofend

\subsection{Minimal uncolorable 
	$3$-uniform bi-hypergraphs}

Now we are going to construct a minimal uncolorable 
$3$-uniform bi-hypergraph of order $n$ for any $n\ge 6$, according to the six cases of $n$ modulo $6$.

\lemm{lem2}
{
For any $k\ge 1$, 
 the bi-hypergraph $\hyh$ 
 obtained from 
$\hyh_{2k}$
by removing  edge $\{\vtx{2k}{i}:i\in[3]\}$ and adding  edges
$\{\vtx{1}{j},\vtx{2k}{j+1},\vtx{2k}{j+2}\}$
for all $j\in [3]$
is minimal uncolorable.
}
\proof
Assume that $c$ is a proper coloring of $\hyh=(\hyv,\hye)$.
Without loss of generality, we can further assume $c(\vtx{1}{1})=c(\vtx{1}{2})=1$ and $c(\vtx{1}{3})=2$. Then by Observation~\ref{obsv2}, we have $c(\vtx{2k-1}{1})=c(\vtx{2k-1}{2})=1$ and $c(\vtx{2k-1}{3})=2$. 
Since $\{v_{2k,j},v_{2k-1,1},v_{2k-1,2}\}$
and $ \{v_{2k,j}, v_{2k-1,j}, v_{2k-1,3}\}$ 
	are edges in $\hyh$
	for both $j=1,2$,
	we have 
	$c(\vtx{2k}{1})=c(\vtx{2k}{2})=2$.
	Then, all vertices in the edge
	$\{v_{1,3}, v_{2k,1}, v_{2k,2}\}$ are assigned color $2$ by $c$, 
	a contradiction
to the assumption of $c$. 
	 Hence $\hyh$ is uncolorable.
	 
Next we show that for any $e\in \hye$, $\hyh-e$ is colorable.
Due to symmetry, we need only to prove 
this conclusion
when
$e=\{\vtx{1}{1},\vtx{2k}{2},\vtx{2k}{3}\}$ or
$$
e\in \bigcup_{i=1}^{2k-1}
\big\{ 
\{ \vtx{i+1}{1},\vtx{i}{1},\vtx{i}{2}\},
\{ \vtx{i+1}{1},\vtx{i}{1},\vtx{i}{3}\},
\{ \vtx{i}{1},\vtx{i}{2},\vtx{i}{3}\}
\big \}.
$$

For each of the four cases, we define a mapping $c$ from $\hyv$ to $\N$ in the following, which can be verified as a proper coloring of $\hyh-e$.

\incase{$e=\{\vtx{1}{1},\vtx{2k}{2},\vtx{2k}{3}\}$.}

For all odd $i\in[2k-1]$ and even $j\in [2k]$,
let $c(\vtx{i}{1})=c(\vtx{j}{2})=c(\vtx{j}{3})=1$ and $c(\vtx{j}{1})=c(\vtx{i}{2})=c(\vtx{i}{3})=2$.

\incase{$e=\{\vtx{q+1}{1},
	\vtx{q}{1},\vtx{q}{2}\}$, where $1\le q\le 2k-1$.}

If $q$ is odd, 
then for all odd $i$ and even $j$
in $[2k]$, let
$$
\left \{
\begin{array}{ll} 
c(\vtx{i}{1})=c(\vtx{i}{2})=c(\vtx{j}{3})=1,
~c(\vtx{j}{1})=c(\vtx{j}{2})=c(\vtx{i}{3})=2,
&\mbox{if }i\le q, j\le q; \\
c(\vtx{i}{1})=c(\vtx{j}{2})=c(\vtx{j}{3})=2,~c(\vtx{j}{1})=c(\vtx{i}{2})=c(\vtx{i}{3})=1,
&\mbox{if }i> q, j>q.
\end{array}
\right. 
$$

Otherwise, $q$ is even. 
Let $c(\vtx{2k}{1})=c(\vtx{2k}{2})=c(\vtx{2k}{3})=3$
and 
for all odd $i$ and even $j$ in $[2k-1]$, let
$$
\left \{
\begin{array}{ll} 
	c(\vtx{i}{1})=c(\vtx{j}{2})=c(\vtx{i}{3})=1,~c(\vtx{j}{1})=c(\vtx{i}{2})
	=c(\vtx{j}{3})=2,
	&\mbox{if }i\le q, j\le q; \\
	c(\vtx{i}{1})=c(\vtx{j}{2})=c(\vtx{j}{3})=3,~c(\vtx{j}{1})=c(\vtx{i}{2})
	=c(\vtx{i}{3})=1,
	&\mbox{if }i> q, j>q.
\end{array}
\right. 
$$

\incase{$e=\{\vtx{q+1}{1},
	\vtx{q}{1},\vtx{q}{3}\}$, where $1\le q\le 2k-1$.}

For all odd $i$ and even $j$ in $[2k]$, let
$$
\left \{
\begin{array}{ll} 
c(\vtx{i}{1})=c(\vtx{j}{2})=c(\vtx{i}{3})=1,~c(\vtx{j}{1})=c(\vtx{i}{2})=c(\vtx{j}{3})=2,
	&\mbox{if }i\le q, j\le q; \\
c(\vtx{i}{1})=c(\vtx{j}{2})=c(\vtx{j}{3})=2,~c(\vtx{j}{1})=c(\vtx{i}{2})=c(\vtx{i}{3})=1,
	&\mbox{if }i> q, j>q.
\end{array}
\right. 
$$

\incase{$e=\{\vtx{q}{1},\vtx{q}{2},
	\vtx{q}{3}\}$, where $1\le q\le 2k-1$.}

If $q$ is odd, 
then let $c(\vtx{q}{1})=c(\vtx{q}{2})=c(\vtx{q}{3})=1$ and
$c(\vtx{2k}{1})=c(\vtx{2k}{2})=c(\vtx{2k}{3})=3$, and 
for all odd $i$ and even $j$ in 
$[2k-1]\setminus \{q\}$, let
$$
\left \{
\begin{array}{ll} 
c(\vtx{i}{1})=c(\vtx{i}{2})=c(\vtx{j}{3})=1,~c(\vtx{j}{1})=c(\vtx{j}{2})=c(\vtx{i}{3})=2, &\mbox{if }i<q, j<q; \\
c(\vtx{i}{1})=c(\vtx{j}{2})=c(\vtx{j}{3})=3,~c(\vtx{j}{1})=c(\vtx{i}{2})=c(\vtx{i}{3})=2,
	&\mbox{if }i>q, j>q.
\end{array}
\right. 
$$

Otherwise, $q$ is even. 
Let $c(\vtx{q}{1})=c(\vtx{q}{2})=c(\vtx{q}{3})=2$, $c(\vtx{2k}{1})=c(\vtx{2k}{2})
=c(\vtx{2k}{3})=3$ and 
for all odd $i$ and even $j$ in $[2k-1]\setminus \{q\}$, let
$$
\left \{
\begin{array}{ll} 
c(\vtx{i}{1})=c(\vtx{i}{2})=c(\vtx{j}{3})=1,~c(\vtx{j}{1})=c(\vtx{j}{2})=c(\vtx{i}{3})=2, &\mbox{if }i<q, j<q; \\
c(\vtx{i}{1})=c(\vtx{j}{2})=c(\vtx{j}{3})=3,~c(\vtx{j}{1})=c(\vtx{i}{2})=c(\vtx{i}{3})=1,
&\mbox{if }i>q, j>q.
\end{array}
\right. 
$$

The result is proven.
\proofend


\lemm{lem1}
{
For any $k\ge 1$,
the bi-hypergraph $\hyh$ 
obtained from $\hyh_{2k+1}$
by removing edge $\{\vtx{2k+1}{i}:
i\in [3]\}$  and adding  edges
$\{\vtx{1}{j},\vtx{2k+1}{j},
	\vtx{2k+1}{j+1}\}$ for all $j\in[3]$
is minimal uncolorable.
}

\proof
Assume that $c$ is a proper coloring of $\hyh=(\hyv,\hye)$, 
 where $c(\vtx{1}{1})=c(\vtx{1}{2})=1$ and $c(\vtx{1}{3})=2$. 
 Then by Observations~\ref{obsv1} and~\ref{obsv2}, we have $c(\vtx{2k}{1})=c(\vtx{2k}{2})=2$ and $c(\vtx{2k}{3})=1$. 
By definition, for each $j\in [2]$, 
$\{\vtx{2k+1}{j}, \vtx{2k}{1}, \vtx{2k}{2}\}$
and $\{\vtx{2k+1}{j}, \vtx{2k}{j}, \vtx{2k}{3}\}$ are edges in $\hyh$, 
implying that $c(\vtx{2k+1}{j})=1$.
It follows that 
$\{\vtx{1}{1},\vtx{2k+1}{1},\vtx{2k+1}{2}\}$ is monochromatic, 
	contradicting the assumption that 
	$c$ is a proper coloring.
Hence $\hyh$ is uncolorable.

Next we show that for any $e\in \hye$, $\hyh-e$ is colorable.
Due to symmetry, we need only to prove this conclusion when $e=\{\vtx{1}{1},\vtx{2k+1}{1},\vtx{2k+1}{2}\}$
or
$$
e\in \bigcup_{i=1}^{2k}
\big \{ \{ \vtx{i+1}{1},\vtx{i}{1},\vtx{i}{2}\},
\{ \vtx{i+1}{1},\vtx{i}{1},\vtx{i}{3}\},
\{ \vtx{i}{1},\vtx{i}{2},\vtx{i}{3}\}
\big \}.
$$

For each of the four cases, we define a mapping $c$ from $\hyv$ to $\N$ in the following, which can be verified as a proper coloring of $\hyh-e$.

\incase{$e=
	\{\vtx{1}{1},\vtx{2k+1}{1},\vtx{2k+1}{2}\}$.}

For all odd $i\in[2k+1]$ and even $j\in [2k]$,
let $c(\vtx{i}{1})=c(\vtx{i}{2})=c(\vtx{j}{3})=1$, $c(\vtx{j}{1})=c(\vtx{j}{2})=c(\vtx{i}{3})=2$.

\incase{$e=\{\vtx{q+1}{1}, 
\vtx{q}{1},\vtx{q}{2}\}$, 
where $1\le q\le 2k$.}

If $q$ is even, then
for all odd $i$ and even $j$ in $[2k+1]$, let
$$
\left \{
\begin{array}{ll} 
c(\vtx{i}{1})=c(\vtx{j}{2})=c(\vtx{i}{3})=1,~c(\vtx{j}{1})=c(\vtx{i}{2})=c(\vtx{j}{3})=2, 
	&\mbox{if }i\le q, j\le q; \\
c(\vtx{i}{1})=c(\vtx{j}{2})=c(\vtx{j}{3})=3,~c(\vtx{j}{1})=c(\vtx{i}{2})=c(\vtx{i}{3})=1,
	&\mbox{if }i>q, j>q.
\end{array}
\right. 
$$

Otherwise, $q$ is odd. 
Let
$c(\vtx{2k+1}{1})=c(\vtx{2k+1}{2})=c(\vtx{2k+1}{3})=3$ and 
	for all odd $i$ and even $j$ in $[2k]$, let
	$$
	\left \{
	\begin{array}{ll} 
		c(\vtx{i}{1})=c(\vtx{j}{2})=c(\vtx{i}{3})=1,~c(\vtx{j}{1})=c(\vtx{i}{2})=c(\vtx{j}{3})=2, 
		&\mbox{if }i\le q, j\le q; \\
		c(\vtx{i}{1})=c(\vtx{j}{2})=c(\vtx{j}{3})=2,~c(\vtx{j}{1})=c(\vtx{i}{2})=c(\vtx{i}{3})=3,
		&\mbox{if }i>q, j>q.
	\end{array}
	\right. 
	$$
	
\incase{$e=\{\vtx{q+1}{1},\vtx{q}{1},\vtx{q}{3}\}$, where $1\le q\le 2k$.}

	For all odd $i$ and even $j$ in $[2k+1]$, let
	$$
	\left \{
	\begin{array}{ll} 
		c(\vtx{i}{1})=c(\vtx{j}{2})=c(\vtx{i}{3})=1,~c(\vtx{j}{1})=c(\vtx{i}{2})=c(\vtx{j}{3})=2, 
		&\mbox{if }i\le q, j\le q; \\
		c(\vtx{i}{1})=c(\vtx{i}{2})=c(\vtx{j}{3})=2,~c(\vtx{j}{1})=c(\vtx{j}{2})=c(\vtx{i}{3})=1,
		&\mbox{if }i>q, j>q.
	\end{array}
	\right. 
	$$
	
\incase{$e=\{\vtx{q}{1}, \vtx{q}{2}, \vtx{q}{3}\}$, where $1\le q\le 2k$.}

If $q$ is odd, then
	let $c(\vtx{q}{1})=c(\vtx{q}{2})=c(\vtx{q}{3})=1$,
		$c(\vtx{2k+1}{1})=c(\vtx{2k+1}{2})=c(\vtx{2k+1}{3})=3$,
	 and 
	for all odd $i$ and even $j$ in $[2k]\setminus \{q\}$, let
	$$
	\left \{
	\begin{array}{ll} 
		c(\vtx{i}{1})=c(\vtx{i}{2})=c(\vtx{j}{3})=1,~c(\vtx{j}{1})=c(\vtx{j}{2})=c(\vtx{i}{3})=2, 
		&\mbox{if }i<q, j<q; \\
		c(\vtx{i}{1})=c(\vtx{i}{2})=c(\vtx{j}{3})=3,~c(\vtx{j}{1})=c(\vtx{j}{2})=c(\vtx{i}{3})=2,
		&\mbox{if }i>q, j>q.
	\end{array}
	\right. 
	$$

Otherwise, $q$ is even. 
	Let $c(\vtx{q}{1})=c(\vtx{q}{2})=c(\vtx{q}{3})=2$,
	$c(\vtx{2k+1}{1})=c(\vtx{2k+1}{2})=c(\vtx{2k+1}{3})=3$,
	and 
	for all odd $i$ and even $j$ in $[2k]\setminus \{q\}$, let
	$$
	\left \{
	\begin{array}{ll} 
		c(\vtx{i}{1})=c(\vtx{i}{2})=c(\vtx{j}{3})=1,~c(\vtx{j}{1})=c(\vtx{j}{2})=c(\vtx{i}{3})=2, 
		&\mbox{if }i<q, j<q; \\
		c(\vtx{i}{1})=c(\vtx{i}{2})=c(\vtx{j}{3})=3,~c(\vtx{j}{1})=c(\vtx{j}{2})=c(\vtx{i}{3})=1,
		&\mbox{if }i>q, j>q.
	\end{array}
	\right. 
	$$

The result is proven.
\proofend

The following four lemmas can be proven in a similar manner. For simplicity, we omit the proofs for minimality.

\lemm{lem3}
{
For any $k\ge 1$,
the bi-hypergraph $\hyh$ 
obtained from $\hyh_{2k}$
by adding a new vertex $v$ and edges 
$\{v,\vtx{1}{j},\vtx{2k}{j+1}\}$ for all $j\in[3]$
is minimal uncolorable.
}
\proof
Assume that $c$ is a proper coloring of $\hyh$. Obviously, there exists a unique $j\in[3]$, such that $c(\vtx{1}{j})$ is unique among $\{c(v_{1,t}):t=1,2,3\}$. By Observations~\ref{obsv1} and~\ref{obsv2}, it can be verified that $c(\vtx{1}{j})=c(\vtx{2k}{j-1})=c(\vtx{2k}{j+1})$ and $c(\vtx{2k}{j})=c(\vtx{1}{j-1})=c(\vtx{1}{j+1})$, where $c(\vtx{1}{j+1})\neq c(\vtx{2k}{j-1})$.
Then, due to the edges $\{v,v_{1,j},v_{2k,j+1}\}$ and $\{v,v_{1,j-1},v_{2k,j}\}$, $c(v)$ must be a value different from $c(\vtx{1}{j})$ and $c(\vtx{1}{j-1})$. However, the edge $\{v,v_{1,j+1},v_{2k,j-1}\}$ indicates that 
$c(v)$ belongs to the set  $\{c(\vtx{1}{j+1}),c(\vtx{2k}{j-1})\}
=\{c(\vtx{1}{j}),c(\vtx{1}{j-1})\}$, a contradiction.
	 Hence $\hyh$ is uncolorable.
\proofend

\lemm{lem4}
{
For any $k\ge 1$,
the bi-hypergraph $\hyh$ 
obtained from $\hyh_{2k+1}$
by adding a new vertex $v$,
removing edge $\{\vtx{1}{i}:i\in [3]\}$ and adding edges $\{v,v_{2k+1,i},v_{2k+1,i+1}\}$
for all $i\in [2]$
and 
$\{v,\vtx{1}{j},
	\vtx{2k+1}{j}\}$ for all $j\in[3]$
is minimal uncolorable.
}
\proof
Assume that $c$ is a proper coloring of $\hyh$. Note that due to edges $\{v_{2,1},v_{2,2},v_{2,3}\}$ and $\{v_{2,j},v_{1,j},v_{1,j+t}\}$ for $j=1,2,3$ and $t=1,2$, $|\{c(v_{1,j}): j=1,2,3\}|$ can only be $1$ or $2$.

Suppose $|\{c(v_{1,j}): j=1,2,3\}|=1$. Then due to the edges $\{v_{2,j},v_{1,j},v_{1,j+t}\}$ for $j=1,2,3$ and $t=1,2$ again, we have $\{c(v_{1,j}): j=1,2,3\}\cap \{c(v_{2,j}): j=1,2,3\}=\emptyset$. Moreover, by Observations~\ref{obsv1} and~\ref{obsv2}, $\{c(v_{1,j}): j=1,2,3\}\cap \{c(v_{2k+1,j}): j=1,2,3\}=\emptyset$. Then, due to $\{v,v_{2k+1,i},v_{2k+1,i+1}\}$ for $i=1,2$, $c(v)$ must belong to $\{c(v_{2k+1,j}): j=1,2,3\}$. Further, there exists $r,l\in[3]$, such that $c(v)=c(v_{2k+1,r})$ and $c(v)\neq c(v_{2k+1,l})$. Since $\left( c(v)\cup c(v_{2k+1,l})\right) \cap \{c(v_{1,j}): j=1,2,3\}=\emptyset$, the edge $\{v,v_{1,l},v_{2k+1,l}\}$ is polychromatic, a contradiction.

Otherwise, $|\{c(v_{1,j}): j=1,2,3\}|=2$. By Observation~\ref{obsv2}, we have $c(v_{1,j})=c(v_{2k+1,j})$ for all $j=1,2,3$. Then, due to the edges $\{v,v_{2k+1,i},v_{2k+1,i+1}\}$ for $i=1,2$, $c(v)=c(v_{2k+1,r})$ for some $r\in [3]$. As a result, $\{v,v_{1,r},v_{2k+1,r}\}$ is monochromatic, a contradiction.

Hence $\hyh$ is uncolorable.
\proofend

\lemm{lem6}
{
	For any $k\ge 1$,
	the bi-hypergraph $\hyh$ 
	obtained from $\hyh_{2k}$
	by adding two new vertices $u$
	and $v$, removing edge $\{\vtx{1}{i}:i\in [3]\}$ and adding edges $\{u,v_{1,1},v_{1,2}\}$, $\{v,v_{1,2},v_{1,3}\}$,
	$\{u,\vtx{1}{j},
		\vtx{2k}{3-j}\}$ and
		$\{v,\vtx{1}{j+1},
		\vtx{2k}{4-j}\}$ for all $j\in[2]$
	is minimal uncolorable.
}

\proof
Assume that $c$ is a proper coloring of $\hyh$. Similar to the proof of Lemma~\ref{lem4}, $|\{c(v_{1,j}): j=1,2,3\}|$ can only be $1$ or $2$.

Suppose $|\{c(v_{1,j}): j=1,2,3\}|=1$. Similarly, we have $\{c(v_{1,j}): j=1,2,3\}\cap \{c(v_{2k,j}): j=1,2,3\}=\emptyset$. 
Obviously, either $c(v_{2k,1})\neq c(v_{2k,2})$ or $c(v_{2k,2})\neq c(v_{2k,3})$ holds. Without loss of generality, assume that $c(v_{2k,1})\neq c(v_{2k,2})$. Then, due to edges $\{u, v_{1,1}, v_{2k,2}\}$ and $\{u, v_{1,2}, v_{2k,1}\}$, we have $c(u)=c(v_{1,1})=c(v_{1,2})$. However, this implies that the edge $\{u,v_{1,1},v_{1,2}\}$ is monochromatic, a contradiction.

Otherwise, $|\{c(v_{1,j}): j=1,2,3\}|=2$. Obviously, either $c(v_{1,1})\neq c(v_{1,2})$ or $c(v_{1,2})\neq c(v_{1,3})$ holds. Without loss of generality, assume that $c(v_{1,1})\neq c(v_{1,2})$. 
By Observations~\ref{obsv1} and~\ref{obsv2}, we have $c(v_{1,1})=c(v_{2k,2})$ and $c(v_{1,2})=c(v_{2k,1})$.
Then, due to the edges $\{u,v_{1,1},v_{2k,2}\}$ and $\{u,v_{1,2},v_{2k,1}\}$, $c(v)$ must be a value different from $c(v_{1,1})$ and $c(v_{1,2})$.
However, the edge $\{u,v_{1,1},v_{1,2}\}$ indicates that $c(u)\in \{c(v_{1,1}),c(v_{1,2})\}$, a contradiction.
Hence $\hyh$ is uncolorable.
\proofend

\lemm{lem5}
{
	For any $k\ge 1$,
	the bi-hypergraph $\hyh$ 
	obtained from $\hyh_{2k+1}$
	by adding 
	two new vertices $u$
	and $v$, 
		removing edge $\{\vtx{1}{i}:i\in [3]\}$ and adding edges $\{u,v_{1,1},v_{1,2}\}$, $\{v,v_{1,2},v_{1,3}\}$,
		$\{u,\vtx{1}{j},
		\vtx{2k+1}{j}\}$ and 
		$\{v,\vtx{1}{j+1},
		\vtx{2k+1}{j+1}\}$ for all $j\in[2]$
	is minimal uncolorable.
}
\proof
The proof is basically the same as that of Lemma~\ref{lem6}.
\proofend 

\vspace{0.2 cm}

Theorem~\ref{thm1} can be obtained by applying Lemmas~\ref{lem2}--\ref{lem5}.

\noindent {\it Proof of Theorem~\ref{thm1}.}
Note that $\hyh_k\in \bH(3k,3,7k-6)$.
 
By Lemma~\ref{lem2}, 
there exists $\hyh=(\hyv, \hye)\in \muc(6k,3)$ 
with $|\hye|=14k-4$, where $k\ge 1$. 
It follows that for $n\in \N$ with 
$n\equiv 0 \pmod{6}$, 
$\muc(n,3)\ne \emptyset$ and 
$\low{n}{3}\le 7\cdot n/3-4$.

By Lemma~\ref{lem1}, 
there exists $\hyh=(\hyv, \hye)\in \muc(6k+3,3)$ 
with $|\hye|=14k+3$, where $k\ge 1$. 
It follows that for $n\in \N$ with 
$n\equiv 3\pmod{6}$, 
$\muc(n,3)\ne \emptyset$ and 
$\low{n}{3}\le 7\cdot n/3-4$.

By Lemma~\ref{lem3}, 
there exists $\hyh=(\hyv, \hye)\in \muc(6k+1,3)$ 
with $|\hye|=14k-3$, where $k\ge 1$. 
It follows that for $n\in \N$ with 
$n\equiv 1\pmod{6}$, 
$\muc(n,3)\ne \emptyset$ and 
$\low{n}{3}\le (7 n-16)/3$.

By Lemma~\ref{lem4}, 
there exists $\hyh=(\hyv, \hye)\in \muc(6k+4,3)$ 
with $|\hye|=14k+5$, where $k\ge 1$. 
It follows that for $n\in \N$ with 
$n\equiv 4\pmod{6}$, 
$\muc(n,3)\ne \emptyset$ and 
$\low{n}{3}\le (7 n-13)/3$.

By Lemma~\ref{lem6}, 
there exists $\hyh=(\hyv, \hye)\in \muc(6k+2,3)$ 
with $|\hye|=14k-1$, where $k\ge 1$. 
It follows that for $n\in \N$ with 
$n\equiv 2\pmod{6}$, 
$\muc(n,3)\ne \emptyset$ and 
$\low{n}{3}\le (7 n-17)/3$.

	By Lemma~\ref{lem5}, 
	there exists $\hyh=(\hyv, \hye)\in \muc(6k+5,3)$ 
	with $|\hye|=14k+6$, where $k\ge 1$. 
	It follows that for $n\in \N$ with 
	$n\equiv 5\pmod{6}$, 
	$\muc(n,3)\ne \emptyset$ and 
	$\low{n}{3}\le (7 n-17)/3$.

Hence the result holds.
\proofend

\section{Further study} 

Note that the value of $\low{6}{3}$ 
follows directly from Theorems~\ref{thm2} and~\ref{thm1}.

\corr{cor4-0}
{$\low{6}{3}=10$.
}

Also, Theorem~\ref{thm2} and Lemma~\ref{lem3} together imply that $10\le \low{7}{3}\le 11$. 
Note that for any 
$\hyh=(\hyv, \hye)$ in $\muc(7,3)$ with $|\hye|=10$, 
 it can be shown by applying the idea in the proof of Lemma~\ref{lem17} 
 that 
$d(v)\ge 4$ for each vertex $v\in \hyv$. We guess that there is no 
bi-hypergraph in $\muc(7,3)$ with $10$ edges. Thus, we propose the following 
conjecture.

\begin{con}
	 $\low{7}{3}=11$.
\end{con}

Motivated by Theorem~\ref{thm2}, we wonder if $m(r)=\binom{(r-1)^2+1}{r}$ holds for every $r\ge 4$.

	\begin{con}\label{con5-3}
		For any $r\in \N$ with $r\ge 4$, $m(r)=\binom{(r-1)^2+1}{r}$.
	\end{con}

For the case $r=4$, 
the above conjecture holds if any $4$-uniform bi-hypergraph with less than $210$ edges is colorable. 
By Lemma~\ref{lem16}, 
 to prove Conjecture~\ref{con5-3}, 
	it remains to show that every $4$-uniform bi-hypergraph of order at most 49 and size at most 209 is colorable.

We also wonder if the conclusion of Lemma~\ref{lem16} holds without restrictions to $r, n$ and $m$.

	\begin{con}\label{con5-4}
	For any $r,m,n\in \N$ with $r\ge 3$ and 
	$n\ge (r-1)^2+1$, 
	if all bi-hypergraphs in $\bH(n,r,m)$ are colorable,
	then all bi-hypergraphs in $\bH(n+1,r,m)$ are also colorable.
\end{con}

Note that Conjecture~\ref{con5-3} follows from Conjecture~\ref{con5-4}.

On the other hand, in the uncolorable $3$-uniform bi-hypergraphs we constructed in Section~\ref{sec1}, many edges overlap heavily, which intrigues us to ask the following problem.

\begin{prob}
	Is every linear uniform 
	bi-hypergraph colorable?
\end{prob}

\section*{Acknowledgement}

This research is supported by 
NSFC (No. 12101347 and 12371340), 
NSF of Shandong Province (No. ZR2021QA085)  and 
the Ministry of Education,
Singapore, under its Academic Research Tier 1 (RG19/22). Any opinions,
findings and conclusions or recommendations expressed in this
material are those of the authors and do not reflect the views of the
Ministry of Education, Singapore.
The first author would like to express her gratitude to National Institute of Education and Nanyang Technological University of Singapore for offering her Nanyang Technological University Research Scholarship during her Ph. D. study.

\end{document}